\documentclass[12pt]{article}
\usepackage{amssymb}

\usepackage{amsmath,latexsym,amsfonts}
\usepackage{graphicx}
\usepackage{amscd}

\newtheorem{lemma}{Lemma}[section]
\newtheorem{coro}[lemma]{Corollary}
\newtheorem{prop}[lemma]{Proposition}
\newtheorem{thm}[lemma]{Theorem}
\newtheorem{defn}[lemma]{Definition}
\makeatletter\@addtoreset{equation}{section}
\renewcommand\theequation{\thesection.\@arabic\c@equation}

\begin{document}
\begin{center}
{\LARGE   Transverse conformal Killing forms and  a Gallot-Meyer
Theorem for foliations}

 \renewcommand{\thefootnote}{}
\renewcommand{\thefootnote}{\arabic{footnote}}
\setcounter{footnote}{0}

\vspace{1 cm} {\large Seoung Dal Jung and Ken Richardson}
\end{center}
\vspace{0.5cm}

{\bf Abstract.} We study transverse conformal Killing forms on
foliations and prove a Gallot-Meyer theorem  for foliations.
Moreover, we show that on a foliation with $C$-positive normal
curvature, if there is a closed basic 1-form  $\phi$ such that
$\Delta_B\phi=qC\phi$, then the foliation is transversally
isometric to the quotient of a $q$-sphere.

\section{Introduction}
On Riemannian manifolds, Killing vector fields and conformal
vector fields are very important geometric objects. Killing forms
and conformal Killing forms are generalizations of such objects
and were introduced by K. Yano [\ref{Yano}] and T.
Kashiwada[\ref{Kashiwada},\ref{Kashiwada1}]; they have been
considered by many researchers
[\ref{Moro},\ref{sem},\ref{tach1},\ref{tach}]. Recently, U.
Semmelmann [\ref{sem}] defined the conformal Killing forms by the
kernel of the twistor operator $T$, which is defined on the
exterior algebra, and proved many interesting results
[\ref{Moro},\ref{sem}].

In this paper, we study the analogous problems on foliations. Let
$\mathcal F$ be a transversally oriented Riemannian foliation on a
compact oriented Riemannian manifold $M$ with codimension $q$. A
transversal infinitesimal automorphism of $\mathcal F$ is an
infinitesimal automorphism preserving the leaves.
 A transversal Killing field is a transversal infinitesimal
isometry, i.e., its flow preserves the transverse metric. A
transversal conformal field is a normal field with a flow
preserving the conformal class of the transverse metric. Such
geometric objects give some information about the leaf space
$M/\mathcal F$. There are several known results about transversal
Killing and conformal fields
[\ref{JJ},\ref{Kamber2},\ref{Pak},\ref{Tond},\ref{Tond1},\ref{Yorozu}].
Since the space of transversal infinitesimal automorphisms can be
identified with the space of the basic 1-forms, we can consider
natural generalizations to differential forms. These are called
transverse conformal Killing forms or transverse twistor forms,
which are defined to be basic forms $\phi$ such that for any
vector field $X$ normal to the foliation,
\begin{align*}
\nabla_X\phi-{1\over r+1}i(X)d\phi+{1\over
q-r+1}X^*\wedge\delta_T\phi=0,
\end{align*}
where $r$ is the degree of the form $\phi$, $X^*$ the dual 1-form
of $X$, and $q$ is the codimension of $\mathcal F$. See Section 3
for the definition of $\delta_T$. The transverse conformal Killing
forms with $\delta_T\phi=0$ are called transverse Killing forms.

This paper is organized as follows. In Section 2, we review
well-known facts concerning basic forms and give a generalization
of Meyer's theorem for foliations. In Section 3, we study
transverse conformal Killing forms on foliations and prove some
vanishing theorems. In Section 4, we study a Gallot-Meyer theorem
[\ref{gallot1}] for foliations. Namely, let $(M,g_M,\mathcal F)$
be a compact Riemannian manifold with a foliation $\mathcal F$
$({\rm codim} \mathcal F=q)$ that has $C$-positive normal
curvature and a bundle-like metric $g_M$. Then, for any basic
$r$-form $\phi$ $( 1\leq r\leq q-1)$, any eigenvalue $\lambda_B$
of the basic Laplacian $\Delta_B$ for $\phi$ satisfies
\begin{equation*}
\lambda_B\geq
 \left\{ \begin{split} &r(q-r+1)C,\quad{\rm if}\ d_B\phi=0,\\
 &(r+1)(q-r)C,\quad{\rm if}\
 \delta_B\phi=0.
 \end{split}
 \right.
 \end{equation*}
Moreover, we show that, if $M$ admits a closed basic 1-form such
that $\Delta_B\phi=qC\phi$ and the foliation $\mathcal F$ has
$C$-positive normal curvature, then the foliation is transversally
isometric to the quotient of a $q$-sphere. Lastly, we study
special transverse Killing forms in Section 5.
 \section{Basic forms and Meyer's Theorem}
 Let $(M,g_M,\mathcal F)$ be a $(p+q)$-dimensional
Riemannian manifold with a foliation $\mathcal F$ of codimension
$q$ and a bundle-like metric $g_M$ with respect to $\mathcal F$.
 Then we have an exact sequence of vector bundles
\begin{align}\label{eq1-1}
 0 \longrightarrow L \longrightarrow
TM {\overset\pi\longrightarrow} Q \longrightarrow 0,
\end{align}
where $L$ is the tangent bundle  and $Q=TM/L$ is the normal bundle
of $\mathcal F$. The metric $g_M$ determines an orthogonal
decomposition $TM=L\oplus L^\perp$, identifying $Q$ with $L^\perp$
and inducing a metric $g_Q$ on $Q$. The metric is bundle-like if
and only if $L_Xg_Q=0$ for every $X\in\Gamma L$, where $L_X$ is
the transverse Lie derivative. Let $V(\mathcal F)$ be the space of
all vector fields $Y$ on $M$ satisfying $[Y,Z]\in\Gamma L$ for all
$Z\in\Gamma L$. An element of $V(\mathcal F)$ is called an {\it
infinitesimal automorphism} of $\mathcal F$. Let
\begin{align}\label{eq1-2}
\bar V(\mathcal F)=\{\bar Y:=\pi(Y)\ |\ Y\in V(\mathcal F)\}.
\end{align}
Then we have an associated exact sequence of Lie algebras
\begin{align}\label{eq1-3}
0\longrightarrow \Gamma L\longrightarrow V(\mathcal
F){\overset\pi\longrightarrow} \bar V(\mathcal F)\longrightarrow
0.
\end{align}
Let $\nabla $ be the transverse Levi-Civita connection on $Q$,
which is torsion-free and metric with respect to $g_Q$. Let
$R^\nabla, K^\nabla,\rho^\nabla$ and $\sigma^\nabla$ be the
transversal curvature tensor, transversal sectional curvature,
transversal Ricci operator and transversal scalar curvature with
respect to $\nabla$, respectively. Let $\Omega_B^*(\mathcal F)$ be
the space of all {\it basic forms} on $M$, i.e.,
\begin{align}\label{eq1-4}
\Omega_B^*(\mathcal F)=\{\omega\in\Omega^*(M)\ | \ i(X)\omega=0,\
i(X)d\omega=0,
                     \quad \forall  X\in \Gamma L\}.
\end{align}
Then $L^2(\Omega^*(M))$ is decomposed as [\ref{Lop}, \ref{Park}]
\begin{align}\label{eq1-5}
L^2(\Omega(M))=L^2(\Omega_B(\mathcal F)) \oplus
L^2(\Omega_B(\mathcal F))^\perp.
\end{align}
We have $\Omega_B^r(\mathcal F)\subset \Gamma(\Lambda^r Q^*)$ and
$\bar V(\mathcal F)\cong \Omega_B^1(\mathcal F)$. Now we define
the connection $\nabla $ on $\Omega _{B}^{\ast }(\mathcal F)$,
which is induced from the connection $\nabla $ on $Q $ and
Riemannian connection $\nabla ^{M}$ of $g_{M}$. This
connection $\nabla $ extends the partial Bott connection $\overset{\circ }{%
\nabla }$ given by $\overset{\circ }{\nabla }_{X}\phi =L_X\phi $
for any $X\in \Gamma L$ [\ref{Kamber1}]. Then the basic forms are
characterized by $\Omega _{B}^{\ast }(\mathcal F)={\rm
Ker}\overset{\circ }{\nabla }\subset \Gamma (\wedge Q^{\ast
}(\mathcal F))$.  Let $P:L^2(\Omega^*(M))\to
L^2(\Omega_B^*(\mathcal F))$ be the orthogonal projection onto
basic forms [\ref{Park}], which preserves smoothness in the case
of Riemannian foliations. For any $r$-form $\phi$, we put the
basic part of $\phi$ as $\phi_B :=P\phi$.
 The exterior differential on the  de Rham complex
$\Omega^*(M)$ restricts a differential $d_B:\Omega_B^r(\mathcal
F)\to \Omega_B^{r+1}(\mathcal F)$.  Let $\kappa\in Q^*$ be the
mean curvature form of $\mathcal F$. Then it is well known
[\ref{Lop}] that $\kappa_B:=P\kappa$ is closed.
 We now recall the star operator $\bar
*:\Omega^r(M)\to \Omega^{q-r}(M)$ given by
[\ref{Park},\ref{Tond1}]
\begin{align}\label{eq1-6}
 \bar * \phi =(-1)^{p(q-r)}*(\phi\wedge \chi_{\mathcal
F}),\quad\forall \phi\in\Omega^r(M),
\end{align}
where $\chi_{\mathcal F}$ is the characteristic form of $\mathcal
F$ and $*$ is the Hodge star operator associated to $g_M$. The
operator $\bar *$ maps basic forms to basic forms and  has the
property that [\ref{Tond1}]
\begin{align}\label{eq1-7}
 * \phi = \bar * \phi \wedge \chi_{\mathcal F},\quad\forall\phi\in\Omega_B^r(\mathcal F).
\end{align}
For any $\phi,\psi\in\Omega_B^r(\mathcal F)$,  $\phi\wedge\bar
*\psi =\psi\wedge\bar *\phi$ and also ${\bar
*}^2\phi=(-1)^{r(q-r)}\phi$. Let $\nu$ be the transversal volume
form, i.e., $*\nu=\chi_\mathcal F$.  The pointwise inner product
$\langle\ , \ \rangle$ on $\Lambda^r Q^*$ is defined uniquely by
\begin{align}\label{eq1-8}
\langle\phi,\psi\rangle\nu=\phi\wedge \bar *\psi.
\end{align}
The global inner product $\ll \cdot,\cdot\gg_B$ on
$L^2(\Omega_B^r(\mathcal F))$ is
\begin{align}\label{eq1-9}
\ll \phi,\psi\gg_B = \int_M \langle\phi,\psi\rangle\mu_M,\quad
\forall\phi,\psi\in\Omega_B^r(\mathcal F),
\end{align}
where $\mu_M=\nu\wedge\chi_{\mathcal F}$ is the volume form with
respect to $g_M$.
 With
respect to this scalar product, the formal adjoint
$\delta_B:\Omega_B^r(\mathcal F)\to \Omega_B^{r-1}(\mathcal F)$ of
$d_B$ is given by [\ref{Park}]
\begin{align}\label{eq1-10}
 \delta_B\phi=(-1)^{q(r+1)+1}\bar *d_T\bar *
 \phi=\delta_T\phi+i(\kappa_B^\sharp)\phi,
 \end{align}
 where $d_T=d-\kappa_B\wedge$ and $\delta_T=(-1)^{q(r+1)+1}\bar
 *d\bar *$ is the formal adjoint operator  of $d_T$ with respect to $L^2(\Omega_B^r(\mathcal F))$. Note that
 \begin{align}\label{eq1-10-1}
 i(\kappa_B^\sharp)\phi=(-1)^{q(r+1)}\bar *\kappa_B\wedge\bar
 *\phi,\quad\forall \phi\in \Omega_B^r(\mathcal F),
 \end{align}
 since $\kappa_B$ is basic.
 \begin{lemma}\label{lemma2-1} On a Riemannian foliation $\mathcal F$, for
 any basic $r$-form $\phi$,
 \begin{align}\label{eq1-100}
 \langle i(\kappa_B^\sharp)\phi,\delta_T\phi\rangle &\geq -\frac12
 |\delta_T\phi|^2-\frac12 |i(\kappa_B^\sharp)\phi|^2,\quad{\rm and}\\
\langle d_B\phi,\kappa_B\wedge\phi\rangle&\leq \frac12|d_B\phi|^2
+\frac12|\kappa_B\wedge\phi|^2.
\end{align}
\end{lemma}
{\bf Proof.} From $|\delta_B\phi|^2\geq 0$ and $|d_T\phi|^2\geq0$,
the proof follows. $\Box$

 The basic
Laplacian $\Delta_B$ is given by
 $ \Delta_B = d_B\delta_B
+ \delta_B d_B$.
 Let $\mathcal H_B^r(\mathcal F)=Ker \Delta_B$
be the set of the {\it basic-harmonic forms} of degree $r$. Then
we have [\ref{Kamber1}, \ref{Park}]
\begin{align}\label{eq1-11}
 \Omega_B^r(\mathcal F) = \mathcal H_B^r(\mathcal F)\oplus{\rm im} d_B \oplus {\rm im} \delta_B
\end{align}
with finite dimensional $\mathcal H_B^r(\mathcal F)$. Let
$\{E_a\}(a=1,\cdots,q)$ be a local orthonormal frame for $Q$,
chosen so that each $E_a\in \bar V(\mathcal F)$. Then the dual
frame $\{\theta^a\}$ consists of local basic forms.  Let
$\nabla_{\rm tr}^*$ be a formal adjoint of $\nabla_{\rm
tr}=\sum_a\theta^a\otimes\nabla_{E_a}:\Omega_B^r(\mathcal F)\to
Q^*\otimes \Omega_B^r(\mathcal F)$. Then $\nabla_{\rm tr}^*
=-\sum_a(i(E_a)\otimes {\rm id})
\nabla_{E_a}+(i(\kappa_B^\sharp)\otimes {\rm id})$, and so
\begin{align}\label{eq1-12}
\nabla_{\rm tr}^*\nabla_{\rm tr} =-\sum_a \nabla^2_{E_a,E_a}
+\nabla_{\kappa_B^\sharp}:\Omega_B^r(\mathcal F)\to
\Omega_B^r(\mathcal F),
\end{align}
where $\nabla^2_{X,Y}=\nabla_X\nabla_Y -\nabla_{\nabla^M_XY}$ for
any $X,Y\in TM$. The operator $\nabla_{\rm tr}^*\nabla_{\rm tr}$
is positive definite and formally self adjoint on the space of
basic forms [\ref{Jung}]. We define the bundle map $A_Y:\Lambda^r
Q^*\to\Lambda^r Q^*$ for any $Y\in V(\mathcal F)$ [\ref{Kamber2}]
by
\begin{align}\label{eq1-13}
A_Y\phi =L_Y\phi-\nabla_Y\phi.
\end{align}
Since $L_X\phi=\nabla_X\phi$ for any $X\in\Gamma L$, $A_Y$
preserves the basic forms and depends only on $\bar Y$. We recall
the generalized Weitzenb\"ock formula.
\begin{thm} $[\ref{Jung1}]$ On a Riemannian foliation $\mathcal F$, we have
\begin{align}\label{eq1-14}
  \Delta_B \phi = \nabla_{\rm tr}^*\nabla_{\rm tr}\phi +
  F(\phi)+A_{\kappa_B^\sharp}\phi,\quad\phi\in\Omega_B^r(\mathcal
  F),
\end{align}
 where $F(\phi)=\sum_{a,b}\theta^a \wedge i(E_b)R^\nabla(E_b,
 E_a)\phi$. If $\phi$ is a basic 1-form, then $F(\phi)^\sharp
 =\rho^\nabla(\phi^\sharp)$.
\end{thm}
Let $\mathcal R^\nabla :\Lambda^2 Q^* \to \Lambda^2 Q^*$ be the
{\it normal curvature operator}, which is defined by
\begin{align}\label{eq1-19}
\langle\mathcal R^\nabla
(\omega_1\wedge\omega_2),\omega_3\wedge\omega_4\rangle=g_Q(R^\nabla(\omega_1^\sharp,\omega_2^\sharp)\omega_4^\sharp,\omega_3^\sharp),
\end{align}
where $\omega_i\in Q^*(i=1,\cdots,4)$.  Then $\mathcal F$ has
constant transversal sectional curvature $C$ if and only if
$\mathcal R^\nabla\omega=C\omega$ for any basic 2-form $\omega$.
\begin{defn}
{\rm A Riemannian foliation $\mathcal F$ is said to have}
$C$-positive normal curvature {\rm if there exists a positive
constant $C$ such that
\begin{align}\label{eq1-20}
\langle\mathcal R^\nabla \omega,\omega\rangle\geq C |\omega|^2
\end{align}
for any 2-form $\omega\in\Gamma (\Lambda^2 Q^*)$.}
\end{defn}
  Then we have a generalization of the Meyer theorem for foliations.
\begin{thm}\label{thm2-5}  Let $\mathcal F$ be a Riemannian
foliation that has $C$-positive normal curvature on a compact
Riemannian manifold $(M,g_M)$. Then for any arbitrary basic
$r$-form $\phi$ $(1\leq r\leq q-1,\ q=codim\mathcal F)$,
\begin{align}\label{eq1-22}
\langle F(\phi),\phi\rangle\geq r(q-r)C|\phi|^2.
\end{align}
The equality holds locally if and only if $\mathcal F$ has
constant transversal sectional curvature $C$.
\end{thm}
{\bf Proof.} The proof is similar to that in [\ref{meyer}]. See
also [\ref{min}].  $\Box$

 Let $\mathcal A$ be the O'Neill's integrability tensor
[\ref{oneil}], which satisfies
\begin{align}\label{eq1-22-1}
\mathcal A_X Y&=\frac12\pi^\perp[X,Y],\quad X,Y\in\Gamma L^\perp,
\end{align}
 where $\pi^\perp:TM\to L$ is an orthogonal projection. It is trivial that $L^\perp$ is integrable if and only if $\mathcal A=0$ on $\Gamma L^\perp$.
Moreover, it is well-known[\ref{oneil}] that, for any unit normal
vectors $X,Y\in\Gamma Q$, we have
\begin{align}\label{eq1-23}
K^\nabla(X,Y)=K^M(X,Y)+3|\mathcal A_XY|^2,
\end{align}
where $K^M$ is the sectional curvature of $g_M$ on $M$.
Equivalently, for any basic 1-forms $\omega_1$ and $\omega_2$,
\begin{align}\label{eq1-23-1}
\langle \mathcal R^\nabla(\omega_1\wedge \omega_2),\omega_1\wedge
\omega_2\rangle =\langle \mathcal R^M(\omega_1\wedge
\omega_2),\omega_1\wedge \omega_2\rangle + |\mathcal
A_{\omega_1^\sharp}{\omega_2^\sharp}|^2,
\end{align}
where $\mathcal R^M$ is the curvature operator on $M$.
 From the equation (\ref{eq1-23-1}), if
$M$ has a $C$-positive curvature, then the foliation $\mathcal F$
has also a $C$-positive normal curvature. Hence we have the
following corollary.
\begin{coro} Let $\mathcal F$ be a Riemannian foliation on a space
$(M,g_M)$ that has $C$-positive curvature and $g_M$ a bundle-like
metric. Then for any basic $r$-form $\phi$,
\begin{align}\label{eq1-24} \langle
F(\phi),\phi\rangle\geq r(q-r)C |\phi|^2.
\end{align}
If the equality holds, then $L^\perp$ is integrable and $M$ has
constant curvature $C$. The converse holds if $M$ has a constant
curvature $C$.
\end{coro}
{\bf Proof.} Inequality (\ref{eq1-24}) is a consequence of Theorem
\ref{thm2-5}, because $\mathcal F$ has $C$-positive normal
curvature. If  equality holds, then the transversal sectional
curvature is the constant $C$. From (\ref{eq1-23}), we know that
$\mathcal A_XY=0$ for all $X,Y\in\Gamma Q$, which means $L^\perp$
is integrable. $\Box$

\section{Transverse conformal Killing forms}
For details for the non-foliation case, see [\ref{sem}].  Let
$(M,g_M,\mathcal F)$ be a Riemannian manifold with a foliation
$\mathcal F$ of codimension $q$ and a bundle-like metric $g_M$.
Then $Q^*\otimes \Lambda^r Q^*$ of $O(q)$-representation is
isomorphic to the following direct sum:
\begin{align}\label{eq2-1}
Q^*\otimes \Lambda^r Q^*\cong \Lambda^{r-1}Q^* \oplus
\Lambda^{r+1}Q^* \oplus \Lambda^{r,1}Q^*,
\end{align}
where $\Lambda^{r,1}Q^*$ is the intersection of the kernels of
wedge product and contraction map. Elements of $\Lambda^{r,1}Q^*$
can be considered as 1-forms on $Q$ with values in $\Lambda^r
Q^*$. For any $s\in Q, \omega\in Q^*$ and $\phi\in\Lambda^r Q^*$,
the projection $Pr_{\Lambda^{r,1}}:Q^*\otimes \Lambda^r Q^*\to
\Lambda^{r,1}Q^*\subset Q^*\otimes \Lambda^r Q^*$ is then
explicitly given by
\begin{align}\label{eq2-2}
[Pr_{\Lambda^{r,1}}(\omega\otimes\phi)](s)=\omega(s)\phi-{1\over
r+1}i(s)(\omega\wedge\phi)-{1\over q-r+1}s^*\wedge
i(\omega^\sharp)\phi,
\end{align}
where $s^*$ denotes the $g_Q$-dual 1-form to $s$, i.e.,
$s^*(v)=g_Q(s,v)$ and $\omega^\sharp$ is the vector defined by
$\omega(s)=g_Q(\omega^\sharp,s)$. Now, we define the operator
$T_{\rm tr}:\Gamma(\Lambda^r Q^*)\to\Gamma( \Lambda^{r,1} Q^*)$ by
\begin{align}\label{eq2-4}
[T_{\rm tr}\phi](v):&=[Pr_{\Lambda^{r,1}}(\nabla_{\rm tr}\phi)](v)\\
&=\nabla_v\phi -{1\over r+1}i(v)d_B\phi +{1\over q-r+1}v^*\wedge
\delta_T\phi\notag
\end{align}
for any $v\in Q$.  The formal adjoint operator $T_{\rm
tr}^*:\Gamma(\Lambda^{r,1}Q^*)\to \Gamma(\Lambda^r Q^*)$ of
$T_{\rm tr}$ is given by
\begin{align}\label{eq2-4-1}
T_{\rm tr}^* = (-\sum_ai(E_a)\otimes id)\nabla_{E_a}
+(i(\kappa_B^\sharp)\otimes id).
\end{align}
For any $\phi\in \Lambda^{r,1}Q^*$ and any $v\in Q$, $v^*\wedge
(T_{\rm tr}\phi)(v)=0=i(v)(T_{\rm tr}\phi)(v)$. So (\ref{eq2-4-1})
is proved.
 \begin{defn} {\rm A basic $r$-form $\phi\in \Omega_B^r(\mathcal F)$ is called a}  transverse conformal
Killing $r$-form {\rm if for any vector field $X\in\Gamma Q$,
\begin{align}\label{eq2-5}
\nabla_X\phi ={1\over r+1}i(X)d_B\phi -{1\over q-r+1}X^*\wedge
\delta_T\phi,
\end{align}
where $\delta_T=\delta_B-i(\kappa_B^\sharp)$. In addition, if the
basic $r$-form $\phi$ satisfies $\delta_T\phi=0$, it is called a}
transverse Killing $r$-form.
\end{defn}
Note that a transverse conformal Killing 1-form is a $g_Q$-dual form of a transverse conformal
Killing vector field. In fact, let $\phi$ be  a transverse conformal Killing 1-form. Then for any
$Y,Z\in \Gamma Q$
\begin{align*}
g_Q(\nabla_Z\phi^\sharp,Y)=(\nabla_Z\phi)(Y) =&\frac12 (i(Z)d_B\phi)(Y) - \frac1q (\delta_T\phi) Z^*(Y)\\
              =&\frac12 d_B\phi(Z,Y) -\frac1q (\delta_T\phi) g_Q(Z,Y).
\end{align*}
Since $\delta_T\phi =-{\rm div}_\nabla \phi^\sharp$, we have
 \begin{equation*}
 (L_{\phi^\sharp}g_Q)(Y,Z)= {2\over q}{\rm div}_\nabla(\phi^\sharp) g_Q(Y,Z)\quad\forall Y,
Z\in\Gamma Q.
\end{equation*}
This means that $\phi^\sharp$ is a transverse conformal Killing
vector field. Note that  since $\phi$ is basic, the $g_Q$-dual
vector $\phi^\sharp$ and $g_M$-dual vector $\phi^\sharp$ are the
same.
\begin{prop}
A basic $r$-form $\phi$ is a transverse conformal Killing form if
and only if, for any vector fields $X,Y, X_a\in \Gamma Q$,
\begin{align*}
&\{i(X)\nabla_Y\phi +i(Y)\nabla_X\phi\}(X_2,\cdots,X_r)\\
&= 2g_Q(X,Y)\theta(X_2,\cdots,X_r)-\sum_{a=2}^r(-1)^a\{
g_Q(Y,X_a)\theta(X,X_2,\cdots,\hat
X_a,\cdots,X_r)\\&+g_Q(X,X_a)\theta(Y,X_2,\cdots,\hat
X_a,\cdots,X_r)\},
\end{align*}
where $\theta=-{1\over q-r+1}\delta_T\phi$ and $\hat X$ means that
$X$ is deleted.
\end{prop}
 A basic $r$-form $\phi$ is a transverse Killing $r$-form if
and only if for any $X\in\Gamma Q$,
\begin{align}\label{eq2-6}
i(X)\nabla_X\phi  =0.
\end{align}
Then we have the following proposition.
\begin{prop} Let $\phi$ be a transverse Killing $r$-form and let
$\gamma$ be a transversal geodesic, i.e., $\gamma'\in \Gamma Q$
and $\nabla_{\gamma'}\gamma'=0$. Then $i(\gamma')\phi$ is parallel
along the transversal geodesic $\gamma$.
\end{prop}
{\bf Proof.} In fact, for any transverse Killing form
$\phi\in\Omega_B^r(\mathcal F)$, we have
\begin{align*}
\nabla_{\gamma'} i(\gamma')\phi
=i(\nabla_{\gamma'}\gamma')\phi+i(\gamma')\nabla_{\gamma'}\phi=0.\quad\Box
\end{align*}
 \begin{lemma} On a Riemannian foliation $\mathcal F$ of codimension $q$, any basic $r$-form
 $\phi$ satisfies
\begin{align}\label{eq2-7}
 |\nabla_{\rm tr}\phi|_B^2 \geq {1\over r+1}|d_B\phi|_B^2 + {1\over
 q-r+1}|\delta_T\phi|_B^2.
 \end{align}
The equality holds if and only if $\phi$ is a transverse conformal
Killing $r$-form, i.e., $T_{\rm tr}\phi=0$.
\end{lemma}
{\bf Proof.} This proof is similar to the one in [\ref {sem}].
$\Box$

 Let $\bar g_M =g_L
+ \bar g_Q (=e^{2u}g_Q)$, where $u$ is a basic function. Then we
have the following.
\begin{prop} On a Riemannian foliation $\mathcal F$, the transverse Levi-Civita connections $\bar
\nabla$ and $\nabla$ are related by
\begin{align}\label{eq3-1}
\bar\nabla_X \phi=\nabla_X\phi -r X(u)\phi -d_B u\wedge i(X)\phi
+X^* \wedge i({\rm grad}_\nabla (u))\phi
\end{align}
 for any $X\in\Gamma Q$ and
$\phi\in\Omega_B^r(\mathcal F)$, where ${\rm grad}_\nabla
(u)=\sum_a E_a(u)E_a$ is a transversal gradient of $u$.
\end{prop}
{\bf Proof.} From the relation between $\nabla$ and $\bar\nabla$
given by
\begin{align*}
\bar\nabla_X\pi(Y) =\nabla_X\pi(Y) +X(u)\pi(Y)
+Y(u)\pi(X)-g_Q(\pi(X),\pi(Y)){\rm grad}_\nabla(u),
\end{align*}
the proof follows. $\Box$

Then we have the following theorem.
\begin{thm} Let $(M,g_M,\mathcal F)$ be a Riemannian manifold
with a foliation $\mathcal F$ and a bundle-like metric $g_M$. Let
$\phi\in\Omega_B^r(\mathcal F)$ be a transverse conformal Killing
$r$-form. Then $\bar\phi=e^{(r+1)u}\phi$ is a transverse conformal
Killing $r$-form with respect to the metric $\bar g_Q=e^{2u}g_Q$.
\end{thm}
{\bf Proof.} First, from (\ref{eq3-1}) and $\kappa_{\bar
g}=e^{-2u}\kappa$ we have that, for any basic $r$-form $\phi$,
\begin{align}\label{eq3-2}
\bar d_B\phi =d_B\phi,\quad \bar \delta_B\phi = e^{-2u}\{\delta_B
+(2r-q)i({\rm grad}_\nabla (u))\}\phi.
\end{align}
 Let $\phi\in\Omega_B^r(\mathcal F)$ be a transverse
conformal Killing $r$-form. Then, for $\bar\phi=e^{(r+1)u}\phi$,
we have from (\ref{eq3-1})
\begin{align*}
i(X)\bar d_B\bar\phi&=(r+1)e^{(r+1)u}\{X(u)\phi-d_Bu\wedge i(X)\phi +{1\over r+1}i(X)d_B\phi\},\\
\hat X^*\wedge\bar
\delta_T\bar\phi&=e^{(r+1)u}X^*\wedge\{\delta_T\phi -(q-r+1)i({\rm
grad}_\nabla (u))\phi\},
\end{align*}
where $\hat X^*=e^{2u}X$ is $\bar g_Q$-dual form to $X$.  Hence
from (\ref{eq2-1}) and (\ref{eq3-1}), we have
\begin{align*}
&{1\over r+1}i(X)\bar d_B\bar\phi -{1\over
q-r+1}\hat X^*\wedge\bar\delta_T\bar\phi\\
&=e^{(r+1)u}\{{1\over r+1}i(X)d_B\phi -{1\over q- r+1}X^*\wedge
\delta_T\phi +X(u)\phi-d_Bu\wedge i(X)\phi\notag\\
 &+X^*\wedge
i({\rm grad}_\nabla (u))\phi\}\\
&= \bar\nabla_X \bar\phi.
\end{align*}
So $\bar\phi$ is a transverse conformal Killing $r$-form with
respect to $\bar g_Q$. $\Box$

 From (\ref{eq2-4-1}),  we have the another generalized
Weitzenb\"ock formula.
 \begin{prop} On a Riemannian foliation $\mathcal F$, we have the
 following Weitzenb\"ock formula;   for any basic $r$-form $\phi$
 \begin{align}\label{eq2-11}
 T_{\rm tr}^*T_{\rm tr}\phi=&\nabla_{tr}^*\nabla_{tr}\phi -{1\over r+1}\delta_B d_B\phi -{1\over
 q-r+1}d_T\delta_T\phi,
 \end{align}
 where $d_T=d_B-\kappa_B\wedge$ is formal adjoint of
$\delta_T$.
\end{prop}
From (\ref{eq1-10}), we have that on $\Omega_B(\mathcal F)$
\begin{align}\label{eq2-11-1}
[T_{\rm tr}^* T_{\rm tr},\bar *]=0.
\end{align}
Hence we have the following corollary.
\begin{coro} Any basic $r$-form $\phi$ is a transverse conformal Killing $r$-form if and only if
$\bar *\phi$ is a transverse conformal Killing $(q-r)$-form.
\end{coro}
 Now we define
$K:\Omega^r_B(\mathcal F)\to \Omega^r_B(\mathcal F)$ by
\begin{align}\label{eq2-9}
K(\phi)=d_B i(\kappa_B^\sharp)\phi +\kappa_B\wedge\delta_T\phi.
\end{align}
Trivially $K$ is formally self adjoint operator and $K$ is
identically zero when $\mathcal F$ is minimal.  From
(\ref{eq1-14}) and (\ref{eq2-11}), we have the following
proposition.
\begin{prop}\label{prop3-9}  On a Riemannian foliation $\mathcal F$, we have that, for any basic $r$-form $\phi\in
\Omega_B^r(\mathcal F)$,
\begin{align}\label{eq2-12}
T_{\rm tr}^*T_{\rm tr}\phi= {r\over r+1}\delta_B d_B\phi + {q-
r\over q-r+1}d_B\delta_B\phi+{1\over q- r+1}K(\phi)
-F(\phi)-A_{\kappa_B^\sharp}\phi.
\end{align}
\end{prop}
\begin{coro} Any basic $r$-form $\phi$
is a transverse conformal Killing form if and only if
\begin{align}\label{eq2-13}
 F(\phi)+A_{\kappa_B^\sharp}\phi
  ={r\over r+1}\delta_B d_B\phi +{q-r\over
  q- r+1}d_B\delta_B\phi+{1\over q- r+1}K(\phi).
 \end{align}
 \end{coro}
 \begin{coro} Any  transverse
 Killing $r$-form $\phi$ satisfies
 \begin{align}\label{eq2-16}
 F(\phi)={r\over r+1}\delta_T d_B\phi.
 \end{align}
\end{coro}
  Now we choose the bundle-like metric $g_M$ such that
 $\delta_B\kappa_B=0$. Any bundle-like metric may be modified to such a metric without changing $g_Q$ [\ref{March},\ref{Mason}]. Then  for any form $\phi$, $\ll\nabla_{\kappa_B^\sharp}\phi,\phi\gg_B=0$.
 Hence any transverse
 conformal Killing $r$-form satisfies
 \begin{align}\label{eq2-14}
 {1\over r+1}\ll i(\kappa_B^\sharp)d_B\phi,\phi\gg_B ={1\over q- r+1}\ll
 \kappa_B\wedge\delta_T\phi,\phi\gg_B.
 \end{align}
  From (\ref{eq2-13}) and (\ref{eq2-14}), any transverse conformal Killing
 $r$-form $\phi$ satisfies
 \begin{align}\label{eq2-17}
&\ll F(\phi),\phi\gg_B\notag\\&={r\over r+1}\Vert d_B\phi\Vert_B^2
+{q- r\over q- r+1}\Vert\delta_T\phi\Vert_B^2+{q-2r\over q-
r+1}\ll \kappa_B\wedge\delta_T\phi,\phi\gg_B.
\end{align}
From (\ref{eq1-100}) and (\ref{eq2-17}), we have
\begin{align}\label{eq2-12-1}
&\ll F(\phi),\phi\gg_B\notag\\&\geq{r\over r+1}\Vert
d_B\phi\Vert_B^2 +{q\over
2(q-r+1)}\Vert\delta_T\phi\Vert_B^2+{2r-q\over 2(q-r+1)}\Vert
i(\kappa_B^\sharp)\phi\Vert_B^2.
\end{align}
 Hence we have the following theorem.
 \begin{thm}\label{prop2-11} Let $(M,g_M,\mathcal F)$ be a compact Riemannian manifold with a foliation $\mathcal F$ of codimension $q$ and a bundle-like metric $g_M$
 such that $\delta_B\kappa_B=0$.
  Suppose $F$ is non-positive and negative at some point. Then, for any $1\leq r \leq q-1$, there are no transverse
conformal Killing $r$-forms on $M$.
\end{thm}
{\bf Proof.} From (\ref{eq2-12-1}), if $2r\geq q$,  it is trivial.
For $2r\leq q$, it follows from the fact that $\bar
*:\Omega_B^r(\mathcal F)\to \Omega_B^{q-r}(\mathcal F)$ is an
isometry and $\bar *\phi$ is a transverse conformal Killing form
when $\phi$ is transverse conformal Killing form. So the proof is
completed. $\Box$
\begin{coro} Let $(M,g_M,\mathcal F)$ be as in Theorem
\ref{prop2-11}. Suppose $F$ is quasi-negative. Then, for any
$1\leq r\leq q-1$, there are no transverse Killing $r$-forms on
$M$.
\end{coro}
Since $F(\phi)^\sharp=\rho^\nabla(\phi^\sharp)$ for any basic
$1$-form $\phi$,  we have the following corollary.
\begin{coro}\label{coro3-3} Let $(M,g_M,\mathcal F)$ be as in Theorem
\ref{prop2-11}. Suppose that $\rho^\nabla$ is non-positive and
negative at some point. Then there are no transverse conformal
Killing vector fields on $M$.
\end{coro}
{\bf Remark.} When $\mathcal F$ is minimal, Corollary
\ref{coro3-3} was proved in [\ref{Kamber2},\ref{Pak}].

\bigskip
\noindent {\bf Remark.} On Riemannian manifolds, the  {\it
conformal Killing forms} are sometimes denoted {\it twistor
forms}, because conformal Killing forms are defined similarly to
the twistor spinors in spin geometry. Moreover, the conformal
Killing form is directly related to the twistor spinor $\psi$,
which satisfies the equation $\nabla^M_X\psi=-\frac1n X\cdot
D\psi$ for vector fields $X$, where $D$ is the Dirac operator. As
with ordinary manifolds, on Riemannian foliations, transverse
conformal Killing forms are related to transversal twistor spinors
$\psi$, which satisfy the equation $\nabla_X\psi=-\frac1q
\pi(X)\cdot D_b\psi -\frac1{2q}\kappa_B\cdot\psi$ for all $X\in
V(\mathcal F)$, where $D_b$ is the basic Dirac operator
[\ref{Jung},\ref{Jung2},\ref{Jung1}]. In fact, let $\psi_1,\psi_2$
be transversal twistor spinors on a transverse spin foliation.
Then the basic $r$-form $\phi_r$ defined on any vectors
$X_1,\cdots,X_r\in V(\mathcal F)$ by
\begin{align}
\phi_r(\bar X_1,\cdots,\bar X_r)=\langle (\bar
X_1^*\wedge\cdots\wedge\bar X_r^*)\cdot\psi_1,\psi_2\rangle
\end{align}
is a transverse conformal Killing form. The proof is similar to
the one in [\ref {sem}].

\section{Gallot-Meyer's theorem for foliations}
Let $(M,g_M,\mathcal F)$ be a  compact Riemannian manifold $M$
with a foliation $\mathcal F$ of codimension $q$ and a bundle-like
metric $ g_M$. We begin with the following Lemma.
\begin{lemma}\label{lemma4-0} For any basic $r$-form $\phi$ on
$M$, we have that
\begin{align}\label{eq4-0}
|\ll \delta_B\phi,i(\kappa_B^\sharp)\phi\gg_B |&\leq
\Vert\delta_B\phi\Vert_B\Vert i(\kappa_B^\sharp)\phi\Vert_B,\\
|\ll d_B\phi,\kappa_B\wedge\phi\gg_B|&\leq\Vert d_B\phi\Vert_B
\Vert\kappa_B\wedge\phi\Vert_B.\label{eq4-0-1}
\end{align}
The equalities hold if and only if $\delta_B\phi=s
i(\kappa_B^\sharp)\phi$ and $d_B\phi=t \kappa_B\wedge\phi$ for
real numbers $s,t$, respectively.
\end{lemma}
{\bf Proof.} The proof follows from the Cauchy-Schwartz
inequality. $\Box$

 From Proposition \ref{prop3-9}, we
have the following proposition.
\begin{prop}\label{prop3-10} Let $(M,g_M,\mathcal F)$ be a compact Riemannian manifold with a foliation $\mathcal F$ of codimension $q$
and a bundle-like metric $g_M$. Then, for any basic $r$-form
$\phi$ $(1\leq r\leq q-1)$,
\begin{align}\label{eq2-19}
&\Vert T_{\rm tr}\phi\Vert_B^2+{1\over \bar r+1}\Vert
i(\kappa_B^\sharp)\phi\Vert_B^2+{\bar r-1\over
\bar r+1}\ll\delta_B\phi,i(\kappa_B^\sharp)\phi\gg_B+\ll i(\kappa_B^\sharp)d_B\phi,\phi\gg_B\notag\\
&={r\over r+1}\Vert d_B\phi\Vert_B^2 +{\bar r\over \bar r
+1}\Vert\delta_B\phi\Vert_B^2-\ll
F(\phi),\phi\gg_B+\frac12\ll(\delta_B\kappa_B)\phi,\phi\gg_B,
\end{align}
where $\bar r = q-r$.
\end{prop}
{\bf Proof.} By direct calculation, $\langle
\nabla_{\kappa_B^\sharp}\phi,\phi\rangle_B=\frac12\kappa_B^\sharp(|\phi|_B^2)=\frac12\langle
d_B|\phi|_B^2,\kappa_B\rangle_B$. By integrating,
$\ll\nabla_{\kappa_B^\sharp}\phi,\phi\gg_B=\frac12\ll(\delta_B\kappa_B)\phi,\phi\gg_B$.
So the proof is completed, using (\ref{eq1-10}), (\ref{eq1-13}),
(\ref{eq2-9}). $\Box$

 From
Proposition \ref{prop3-10}, we can  generalize the Gallot-Meyer's
Theorem(cf. [\ref{gallot1}]) to the foliation setting.
\begin{prop}\label{thm3-14} Let $(M,g_M,\mathcal F)$ be a compact Riemannian manifold with a foliation $\mathcal F$ $(codim \mathcal F=q)$
that has $C$-positive normal curvature and a bundle-like metric
$g_M$. Then, for any basic $r$-form $\phi$ $( 1\leq r\leq q-1)$,
any eigenvalue $\lambda_B$ of $\Delta_B$ for $\phi$ satisfies
\begin{equation}\label{eq4-7}
\lambda_B\geq
 \left\{ \begin{split} &r(q-r+1)C+B_1,\quad{\rm if}\ d_B\phi=0,\\
 &(r+1)(q-r)C+B_2,\quad{\rm if}\
 \delta_B\phi=0,
 \end{split}
 \right.
 \end{equation}
 where
\begin{align*}
 B_1&={\alpha_1^2\over
2}-\alpha_2-\alpha_1\sqrt{r(q-r+1)C +{\alpha_1^2 \over
4}-\alpha_2},\\
B_2&={\beta_1^2 \over 2}-\beta_2 -\beta_1\sqrt{(r+1)(q-r)C
+{\beta_1^2 \over 4}-\beta_2},
\end{align*}
$\alpha_1={q-r-1\over q-r} {\rm max}(|\kappa_B|)$,
$\alpha_2={q-r+1\over 2(q-r)}{\rm max }(\delta_B\kappa_B)$,
$\beta_1={r+1\over r}{\rm max}(|\kappa_B|)$ and $\beta_2={r+1\over
2r}{\rm max}(\delta_B\kappa_B)$.
\end{prop}
{\bf Proof.} Let $\Delta_B\phi=\lambda_B\phi$ and $d_B\phi=0$.
Since $\Vert \kappa_B\wedge\phi\Vert_B^2 +\Vert
i(\kappa_B^\sharp)\phi\Vert_B^2=\Vert |\kappa_B|\phi\Vert_B^2$,
from (\ref{eq4-0}), we have
\begin{align}\label{eq2-20}
|\ll i(\kappa_B^\sharp)\phi,\delta_B\phi\gg_B|\leq
\lambda_B^{1\over 2}\Vert |\kappa_B|\phi\Vert_B
\Vert\phi\Vert_B\leq \lambda_B^{1\over 2}{\rm
max}(|\kappa_B|)\Vert\phi\Vert_B^2.
\end{align}
From (\ref{eq1-22}), (\ref{eq2-19}) and (\ref{eq2-20}), we have
\begin{align}\label{eq2-23}
&\Vert T_{\rm tr}\phi\Vert_B^2 +{1\over q-r+1}\Vert
i(\kappa_B^\sharp)\phi\Vert_B^2\\
&\leq\int_M\Big({q-r\over q-r+1}\lambda_B +{q-r-1\over
q-r+1}|\kappa_B|\lambda_B^{1\over
2}-r(q-r)C+\frac12(\delta_B\kappa_B)\Big)|\phi|^2\notag\\
 &\leq \Big({q-r\over q-r+1}\lambda_B +{q-r-1\over q-r+1}{\rm
max}(|\kappa_B|)\lambda_B^{1\over 2}-r(q-r)C +{1\over 2}{\rm
max}(\delta_B\kappa_B)\Big)\Vert \phi\Vert_B^2\notag\\
&={q-r\over q-r+1}\{\lambda_B +\alpha_1\lambda_B^{1\over
2}-r(q-r+1)C+\alpha_2\}\Vert\phi\Vert_B^2,\notag
\end{align}
where $\alpha_1={q-r-1\over q-r} {\rm max}(|\kappa_B|)\geq 0$ and
$\alpha_2={q-r+1\over 2(q-r)}{\rm max }(\delta_B\kappa_B)\geq 0$.
Hence we have
\begin{align}
\lambda_B +\alpha_1\lambda_B^{1\over 2}-r(q-r+1)C+\alpha_2\geq 0,
\end{align}
which yields
\begin{align}\label{eq6-3}
\lambda_B^{1\over 2}\geq \sqrt{r(q-r+1)C +{\alpha_1^2 \over
4}-\alpha_2}-{\alpha_1\over 2}.
\end{align}
Hence we have
\begin{align}
\lambda_B \geq r(q-r+1)C + {\alpha_1^2\over
2}-\alpha_2-\alpha_1\sqrt{r(q-r+1)C +{\alpha_1^2 \over
4}-\alpha_2}.
\end{align}
Therefore, the first inequality of (\ref{eq4-7}) is proved. For
the proof of the second inequality of (\ref{eq4-7}), let
$\Delta_B\phi=\lambda_B\phi$ and $\delta_B\phi=0$. From
(\ref{eq4-0}), we get
\begin{align}\label{eq6-4-1}
|\ll d_B\phi,\kappa_B\wedge\phi\gg_B|\leq \lambda_B^{1\over 2}{\rm
max}(|\kappa_B|) \Vert \phi\Vert_B^2.
\end{align}
From (\ref{eq1-22}), (\ref{eq2-19}) and (\ref{eq6-4-1}), we have
\begin{align}\label{eq6-4-2}
&\Vert T_{\rm tr}\phi\Vert_B^2 +{1\over q-r+1}\Vert
i(\kappa_B^\sharp)\phi\Vert_B^2\\
 &\leq \Big({r\over r+1}\lambda_B +{\rm
max}(|\kappa_B|)\lambda_B^{1\over 2}-r(q-r)C +{1\over 2}{\rm
max}(\delta_B\kappa_B)\Big)\Vert \phi\Vert_B^2\notag\\
&={r\over r+1}\{\lambda_B +\beta_1\lambda_B^{1\over
2}-(r+1)(q-r)C+\beta_2\}\Vert\phi\Vert_B^2,\notag
\end{align}
where $\beta_1={r+1\over r}{\rm max}(|\kappa_B|)\geq 0$ and
$\beta_2={r+1\over 2r}{\rm max}(\delta_B\kappa_B)\geq 0$. Hence we
have
\begin{align}\label{eq6-4-3}
\lambda_B^{1\over 2} \geq \sqrt{(r+1)(q-r)C +{\beta_1^2 \over
4}-\beta_2}-{\beta_1\over 2},
\end{align}
which yields
\begin{align}
\lambda_B\geq (r+1)(q-r)C +{\beta_1^2 \over 2}-\beta_2
-\beta_1\sqrt{(r+1)(q-r)C +{\beta_1^2 \over 4}-\beta_2}.
\end{align}
Hence (\ref{eq4-7}) follows. $\Box$

\begin{prop} \label{thm4-5} Let $(M,g_M,\mathcal F)$ be as in
Proposition \ref{thm3-14}. Then, for any basic $r$-form $\phi$
$(0<r<q)$,

 $(1)$ if $\phi\in{\rm Ker} d_B$ is an eigenform
corresponding to $\lambda_B=r(q-r+1)C+B_1$, then $\phi$ is a
transverse conformal Killing $r$-form and $\kappa_B=0$.

$(2)$ if $\phi\in  {\rm Ker} \delta_B$ is an eigenform
corresponding to $\lambda_B=(r+1)(q-r)C+B_2$, then $\phi$ is a
transverse Killing $r$-form and $\kappa_B=0$.
\end{prop}
{\bf Proof.} Let $\phi\in {\rm Ker}d_B$ be an eigenform with
$\lambda_B=r(q-r+1)C+B_1$.  From (\ref{eq2-23}), $T_{\rm
tr}\phi=0$ and $i(\kappa_B^\sharp)\phi=0$. So $\phi$ is the
transverse conformal Killing form. Moreover, from the second
inequality in (\ref{eq2-23}), we have
\begin{align*}
|\kappa_B|={\rm max}(|\kappa_B|),\quad \delta_B\kappa_B={\rm
max}(\delta_B\kappa_B).
\end{align*}
So $|\kappa_B|$ and $\delta_B\kappa_B$ are constant. Since
$d_B\kappa_B=0$ [\ref{Lop}], we have $\Delta_B\kappa_B=0$, i.e.,
$\kappa_B\in \mathcal H_B^1(\mathcal F)$. On a foliation $\mathcal
F$ that has $C$-positive normal curvature, $\mathcal
H_B^1(\mathcal F)=0$ [\ref{min}].  So there exists a basic
function $h$ such that $\kappa_B=dh$. Therefore, $0=\int_M\Delta
h=\int_M(\delta_B\kappa_B)=(\delta_B\kappa_B){\rm Vol}(M)$, which
mean $\delta_B\kappa_B=0$. This means that $h$ is a harmonic
function on $M$, and so $h$ is constant, since $M$ is compact.
Hence $\kappa_B=0$. This proves (1). The proof of (2) is similar.
$\Box$

\vskip 0.5cm
 \noindent{\bf Remark.} Note that $B_1(\leq 0)$ and
$B_2(\leq 0)$ depend on the mean curvature form of
 $\mathcal F$. Therefore, Proposition 4.4 implies that, on a foliation
$\mathcal F$ that has $C$-positive normal curvature, if there
exists a closed basic $r$-form (resp. coclosed basic $r$-form)
with eigenvalue $\lambda_B=r(q-r+1)C +B_1$ (resp.
$(r+1)(q-r)C+B_2$), then $B_1=0$ (resp. $B_2=0$).

\bigskip
  Now, we recall the tautness theorem [\ref{Lop},\ref{mmr}] on a foliated Riemannian
manifold.
\begin{thm} \label{thm4-5}$(${\bf Tautness theorem}$)$ Let $(M,g_M,\mathcal F)$ be a compact Riemannian manifold with a
Riemannian foliation $\mathcal F$ of codimension $q\geq 2$ and a
bundle-like metric $g_M$. If the transversal Ricci operator
$\rho^\nabla$ is positive definite, then $\mathcal F$ is taut,
i.e., there exists a bundle-like metric $\tilde g_M$  for which
all leaves are minimal submanifolds.
\end{thm}
From Proposition 4.3 and Theorem \ref{thm4-5},  we have the
following theorem.
\begin{thm} Let $(M,g_M,\mathcal F)$ be as in Theorem
\ref{thm3-14}. Then, for any basic $r$-form $\phi$ $( 1\leq r\leq
q-1)$, any eigenvalue $\lambda_B$ of $\Delta_B$ for $\phi$
satisfies
\begin{equation}\label{eq4-71}
\lambda_B\geq
 \left\{ \begin{split} &r(q-r+1)C,\quad{\rm if}\ d_B\phi=0,\\
 &(r+1)(q-r)C,\quad{\rm if}\
 \delta_B\phi=0.
 \end{split}
 \right.
 \end{equation}
\end{thm}
{\bf Proof.} From the assumption $\mathcal R^\nabla\geq C\cdot
id$, we have $\rho^\nabla\geq (q-1)C \cdot id$. From Theorem 4.5,
 since $\mathcal{F}$ is taut, the basic component of the mean
curvature form $\kappa$ is exact. By the reasoning in the proof of
Corollary 4.2 in [13], we may modify the bundle-like metric
$g_{M}$ such that the basic Laplacian is unchanged as an operator,
such that the transverse metric is the same as that of the
original metric, and such that $\kappa=0$. Hence the proof follows
from Proposition 4.3. $\Box$

\vskip 0.5cm \noindent{\bf Remark.} Let $\mathcal F$ have
$C$-positive normal curvature. Then, by the tautness theorem, we
can choose a bundle-like metric $g_M$ such that $\kappa=0$. So,
from Theorem \ref{thm4-5}, any closed $r$-form $\phi$ with
$\lambda_B\phi=r(q-r+1)C$ is a transverse conformal Killing
$r$-form, and a coclosed $r$-form $\phi$ with
$\lambda_B\phi=(r+1)(q-r)C$ is a transverse Killing $r$-form.

\bigskip

 Now, we recall the generalized Obata
theorem [\ref{ken}] for foliations. Let $\rho^\nabla(X)\geq
(q-1)CX$ for any constant $C>0$. If the smallest eigenvalue
$\lambda_1$ of the basic Laplacian acting on functions is
$\lambda_1=qC$, then $\mathcal F$ is minimal and transversally
isometric to the action of a finite subgroup of $O(q)$ acting on
the $q$-sphere of constant curvature $C$. Hence we have the
following theorem.
\begin{thm} Let $(M,g_M,\mathcal F)$ be a compact Riemannian
manifold with a foliation $\mathcal F$ that has $C$-positive
normal curvature  and a bundle-like metric $g_M$. If $M$ admits a
closed basic 1-form $\phi$ such that $\Delta_B\phi=qC\phi$, then
$\mathcal F$ is minimal and transversally isometric to the action
of a finite subgroup of $O(q)$ acting on the $q$-sphere of
constant curvature $C$.
\end{thm}
{\bf Proof.} Let $\phi$ be a closed basic 1-form and
$\Delta_B\phi=qC\phi$. If we put $f=\delta_B\phi$, then
\begin{align*}
\Delta_B f=\delta_B d_B\delta_B\phi
=\delta_B\Delta_B\phi=qC\delta_B\phi=qC f.
\end{align*}
Since the normal curvature operator satisfies $\mathcal R^\nabla
\geq C\cdot id$ on $\Lambda^2Q^*$, we have $\rho^\nabla(X)\geq
(q-1)C X$. So by the generalized Obata theorem, the foliation is
minimal and transversally isometric to the action of a finite
subgroup of $O(q)$ acting on the $q$-sphere of constant of
curvature $C$. $\Box$

\section{ Special transverse Killing forms}

\begin{defn} A transverse Killing $r$-form $\phi$ is called special with
$\beta$ if it satisfies
\begin{align}\label{eq4-1}
\nabla_X d_B\phi= \beta X^*\wedge\phi,\quad\forall X\in\Gamma (Q),
\end{align}
where $\beta$ is a constant.
\end{defn}
Let $\phi\in\Omega_B^r(\mathcal F)$ be a special transverse
Killing $r$-form with $\beta$. Then we have
\begin{align}\label{eq4-2}
\Delta_B\phi=\delta_Td_B\phi +L_{\kappa_B^\sharp}\phi.
\end{align}
Since $\delta_T=-\sum_a i(E_a)\nabla_{E_a}$  on
$\Omega_B^*(\mathcal F)$, from (\ref{eq4-2}) we have
\begin{align}\label{eq4-3}
\Delta_B\phi-L_{\kappa_B^\sharp}\phi=-\beta(q-r)\phi.
\end{align}
From ({\ref{eq2-16}), (\ref{eq4-2}) and (\ref{eq4-3}), if $\phi$
is a special transverse Killing $r$-form with $\beta$, then
\begin{align}\label{eq4-4}
F(\phi)=-{r(q-r)\over r+1}\beta\phi.
\end{align}
So we have the following proposition.
\begin{prop} Let $\mathcal F$ be a Riemannian foliation that has $C$-positive normal curvature
on a Riemannian manifold $(M,g_M)$. If $\phi$ is a special
transverse Killing $r$-form with $\beta$, then
\begin{align*}
 \beta\leq -(r+1)C.
 \end{align*}
\end{prop}
From (\ref{eq4-4}) and Theorem \ref{thm2-5}, we have the following
theorem.
\begin{thm} Let $(M,g_M,\mathcal F)$ be a compact Riemannian
manifold with a foliation $\mathcal F$ of codimension $q$ and a
bundle-like metric $g_M$. Assume that the transversal sectional
curvature $K^\nabla$ is a positive constant $C$. Then any
transverse Killing $r$-form is special with $\beta=-(r+1)C$.
\end{thm}
{\bf Proof.} Let $\phi$ be a transverse Killing $r$-form, i.e.,
$\nabla_X\phi={1\over r+1} i(X)d_B\phi$. Hence we have
\begin{align*}
R^\nabla(X,Y)\phi={1\over
r+1}\{i(Y)\nabla_X-i(X)\nabla_Y\}d_B\phi.
\end{align*}
Hence we have
\begin{align}\label{eq4-5}
\sum_a\theta^a\wedge R^\nabla(X,E_a)\phi={r\over r+1}\nabla_X
d_B\phi.
\end{align}
Since $\mathcal F$ has constant transversal curvature $C$, we have
from
\begin{align}\label{eq4-6}
\sum_a\theta^a\wedge R^\nabla(X,E_a)\phi=-rC X^*\wedge\phi.
\end{align}
From (\ref{eq4-5}) and (\ref{eq4-6}), the proof is completed.
$\Box$

\bigskip
\noindent{\bf Acknowledgements} This paper was supported by
KRF-2007-313-C00064 from Korea Research Foundation.

\vskip 1.0cm Seoung Dal Jung

\noindent Department of Mathematics, Cheju National University,
Jeju 690-756, Korea(e-mail: sdjung@cheju.ac.kr)

\vskip 0.5cm
\noindent Ken Richardson

\noindent Department of Mathematics, Texas Christian University,
Fort Worth, TX 76129, U.S.A.(e-mail:k.richardson@tcu.edu)

\end{document}